\documentclass{IEEEtran}
\usepackage{amsmath}
\usepackage{amsfonts}
\usepackage{xcolor}
\usepackage{graphicx}

\usepackage{algorithm}
\usepackage{algorithmic}
\usepackage{booktabs}
\usepackage[sort]{cite}
\usepackage[keeplastbox]{flushend}
\newtheorem{theorem}{Theorem}
\newtheorem{assumption}{Assumption}
\newtheorem{lemma}{Lemma}
\newtheorem{remark}{Remark}
\newtheorem{corollary}{Corollary}
\newenvironment{proof}{{\indent\indent\itshape Proof:}\,}{\hfill$\square$}

\newcommand{\normsq}[1]{\|#1\|^2}
\newcommand{\norm}[1]{\|#1\|}

\newcommand{\defeq}{\stackrel{\triangle}{=}}
\newcommand{\paren}[1]{\left(#1\right)}

\newcommand{\Ucal}{\mathcal{U}}

\newcommand{\kth}{k{\text{th}}}

\newcommand{\OMIT}[1]{}
\newcommand*{\MATRIX}[1]{\ensuremath{\mathbf{#1}}}

\definecolor{RED}{rgb}{1,0,0}
\newcommand{\dopr}[1]{\text{d}{#1}}
\newcommand{\secondRevision}[1]{\textcolor{black}{#1}}
\newcommand{\thirdRevision}[1]{\textcolor{black}{#1}}

\newif\ifFullVersion
\FullVersiontrue
\begin{document}

\title{Integral-Type Event-Triggered Model Predictive Control of Nonlinear Systems with Additive Disturbance} % Title, preferably not more 
% than 10 words.
\author{Qi~Sun,~\IEEEmembership{Student Member,~IEEE},
	Jicheng~Chen,~\IEEEmembership{Student Member,~IEEE},
	Yang~Shi,~\IEEEmembership{Fellow,~IEEE}% <-this % stops a space
	\thanks{Corresponding author Y.~Shi. Tel. +1-250-853-3178. Fax. +1-250-721-6051. The authors are with the Department of Mechanical Engineering, University of Victoria, Victoria, BC V8W 2Y2, Canada.}}% <-this % stops a space
\maketitle

\begin{abstract}                          % Abstract of not more than 200 words.
	This paper studies integral-type event-triggered model predictive control (MPC) of continuous-time nonlinear systems. An integral-type event-triggered mechanism is proposed by incorporating the integral of errors between the actual and predicted state sequences, leading to reduced average sampling frequency. Besides, a new and improved robustness constraint is introduced to handle the additive disturbance, rendering the MPC problem with a potentially enlarged initial feasible region. Furthermore, the feasibility of the designed MPC and the stability of the closed-loop system are rigorously investigated. Several sufficient conditions to guarantee these properties are established, which is related to factors such as the prediction horizon, the disturbance bound, the triggering level, and the contraction rate for the robustness constraint. The effectiveness of the proposed algorithm is illustrated by numerical examples and comparisons.
\end{abstract}

                             % keyword list or with the 
% help of the Automatica 
% keyword wizard

\section{Introduction}
Recent research interests and efforts have been directed towards reducing the communication load in control systems by using the event-triggered scheme~\cite{aastrom2002comparison,tabuada2007event,lunze2010state,donkers2012output,heemels2013periodic,mousavi2016integral,yu2017input,hu2015consensus}. Such a control paradigm employs the so-called event-triggered mechanism (ETM) to reduce the frequency of state sampling and information transmission while maintaining expected performance. Compared with the conventional periodic control, event-triggered control treats sampling instants as a design factor, whereas the periodic scheme samples the states at fixed time instants. The event-triggered paradigm features the design of an appropriate scheduling mechanism with the capability of determining when the states or sensor outputs should be sampled. The benefit of such an effective scheduling can lead to a reduced communication rate and thus decrease the communication load, especially in cases where the communication resources are limited. 
\par 
Several pioneering works have been devoted to building the fundamental blocks of the event-triggered control~\cite{aastrom2002comparison,tabuada2007event}. The presented approach~\cite{aastrom2002comparison} used a constant threshold-based event-triggered scheme, which executes the next sampling when the norm of state or estimation error exceeds a certain constant bound. Compared with the conventional periodic control, this method is shown to have great advantages in terms of reducing communication rate~\cite{aastrom2002comparison}. In another early work~\cite{tabuada2007event}, the proposed event-triggered scheme adopted a so-called relative threshold policy for a class of nonlinear systems. In addition, a lower bound for the inter-execution time is proved to exist for avoiding the Zeno behavior. The study of state-based event-triggered control can be found in~\cite{lunze2010state}, and output-based event-triggered control has been reported in~\cite{donkers2012output}. In~\cite{heemels2013periodic}, the authors proposed a periodic event-triggered scheme for reducing the communication rate, where the event-triggering condition is measured periodically with a fixed time interval. To further reduce the communication rate, the authors in~\cite{mousavi2016integral} proposed a novel integral-based event-triggered scheme by incorporating the integral of estimated errors to the event-triggering condition. Input-to-state stability of the integral-based event-triggered control is investigated in~\cite{yu2017input}. \thirdRevision{Event-triggered schemes have also found applications in multi-agent systems, such as consensus~\cite{hu2015consensus,mu2018event}, distributed formation control~\cite{ge2017distributed}, and distributed estimation~\cite{ge2019distributed}.}
\par 
In control applications, model predictive control (MPC) has been one of the most successful control methodologies. The basic and essential idea of MPC framework is to solve optimization problems online at each sampling instant, and apply the optimized control action to the plant. The advantages of MPC is that it can simultaneously handle the system constraints and take the future system behavior into account. \secondRevision{Compared with other well-known networked control strategies such as networked PI control~\cite{zhang2019networkPIcontrol} and event-triggered fuzzy control~\cite{zhang2015networkFuzzyControl}, the MPC-based method provides constraint satisfactions and better performance guarantees in the control design of networked control systems~\cite{ahn2014receding,varutti2009event,eqtami2011event,li2014event}}.
\par
In particular, the incorporation of the event-triggered scheme into MPC can significantly alleviate the communication and computation load, especially considering that MPC features relatively heavier computational complexity. Therefore, event-triggered MPC (ET-MPC) has received many research efforts~\cite{varutti2009event,eqtami2011event,eqtami2011novel,eqtami2012event,lehmann2013event,li2014event,li2015neighbor,hashimoto2016self,brunner2017robust,liu2018robust,li2018triggering}. Specifically, by using the event-triggered scheme, ET-MPC can reduce the frequency of state sampling, optimization solving, and information transmission. In the literature, the research on ET-MPC can be generally classified into two categories, e.g., for linear~\cite{lehmann2013event,brunner2017robust} and nonlinear systems~\cite{eqtami2011novel,li2014event,hashimoto2016self,liu2018robust,li2018triggering}, respectively. (1) {\bf ET-MPC for linear systems.} The authors in~\cite{brunner2017robust} studied the ET-MPC of linear systems with additive disturbance by using a tube-based approach. (2) {\bf ET-MPC for nonlinear systems.} In~\cite{eqtami2011novel}, an event-triggered scheme was proposed for nonlinear systems with additive disturbance by continuously measuring the discrepancy between the actual and the predicted trajectories. In order to acquire the benefit of avoiding the Zeno behavior, the authors in~\cite{li2014event} proposed an event-triggered mechanism design which can guarantee that the inter-execution time is lower bounded. In~\cite{hashimoto2016self}, a self-triggered MPC framework is proposed for nonlinear affine systems to further reduce the communication load. The control signal is a piecewise constant signal by using sample-and-hold on the optimal control sequence, where each sampling interval is chosen by an adaptive selection scheme. It is worthwhile to note that ET-MPC has also been applied to decentralized systems~\cite{eqtami2011event} and distributed systems~\cite{eqtami2012event,li2015neighbor}.
\par 
In this paper, we investigate the {\bf \em integral-type} {\bf ET-MPC} for continuous-time nonlinear systems with additive disturbance, aiming at alleviating the communication load while ensuring the feasibility of proposed MPC algorithm and the stability of the closed-loop system. The main contributions and novelties of this work are three-fold: 1) An integral-type ET-MPC algorithm for the continuous-time nonlinear system with additive disturbance is proposed, where the integral of errors between the actual and predicted state sequences is used as the triggering condition for saving more communication resources. 2) A new and improved robustness constraint on the system states, rendering the MPC problem with a potentially enlarged initial feasible region compared with the existing one, is proposed to the MPC problem. 3) The feasibility of the designed MPC and the stability of the closed-loop system are thoroughly studied, theoretically showing that both of these important properties are related to the factors including the prediction horizon, the bound of disturbance, the triggering level, and the contraction rate of the robustness constraint.
%Moreover, there exists a design trade-off for the parameters such as triggering level $\delta$ and prediction horizon $T$ when the control performance and the communication load are considered.
\par
The remainder of paper is organized as follows. Section~\ref{sec:problemFormulation} describes the problem formulation. Section~\ref{sec:INTET-MPC} proposes the integral-type ET-MPC scheme. Section~\ref{sec:mainresult} presents feasibility of the MPC problem and stability of the closed-loop system. Section~\ref{sec:simulationresult} illustrates a simulation example to verify the effectiveness of the proposed scheme and algorithm. Finally, we conclude this work in Section~\ref{sec:conclusion}.
\par 
\noindent
%\section*{Notions}
\textbf{Notations:} The real space is denoted by $\mathbb{R}$ and the set of all positive integers is given by $\mathbb{N}$. For a given matrix $X$, $X^{\top}$ and $X^{-1}$ denote its transpose and inverse (if invertible), respectively. For a symmetric matrix $S\in\mathbb{R}^{n\times n}$, $S\succ 0$ and $S\succeq 0$ is used as a common notation for positive definite (PD) and positive semidefinite (PSD) matrices; the largest and smallest eigenvalues of $S$ are denoted by $\overline{\lambda}(S)$ and $\underline{\lambda}(S)$. Given a column vector $x\in\mathbb{R}^{n}$, $\norm{x}:=\sqrt{x^{\top}x}$ represents the Euclidean norm and $\norm{x}_{P}:=\sqrt{x^{\top} P x}$ is the $P$-weighted norm. 
% We also use the notation $\normsq{x}_{P}:=x^{\top} P x$. 
% Given a differentiable vector-valued function $g(t)$ on $[a,b]$, we use $g^{'}(t)$ to represent its Jacobian matrix. 

\section{Problem Formulation} \label{sec:problemFormulation}
% \subsection{System Dynamics}
We consider a continuous-time nonlinear system with additive disturbance as follows
\begin{equation} \label{eq:nonlinearDynamics}
\dot{x}(t) = f(x(t),u(t)) + \omega(t), 
\end{equation}
where $x(t)\in \mathbb{R}^n$ is the state variable, $u(t)\in\mathbb{R}^m$ is the control input, and $\omega(t)\in\mathbb{R}^n$ is \thirdRevision{unknown but bounded additive disturbance}. The system satisfies $f(0,0)=0$ and has a Lipschitz constant $L$. The control input constraint is $u(t)\in \mathcal{U}$, where $\mathcal{U}\in\mathbb{R}^m$ is a compact set containing the origin. Moreover, the disturbance $\omega(t)$ belongs to a compact set $\mathcal{W}$ and its upper bound is given by $\rho\defeq \sup_{\omega(t)\in\mathcal{W}}\norm{\omega(t)}$. The nominal system of~(\ref{eq:nonlinearDynamics}) can be given as
\begin{equation} \label{eq:NominalNonlinearDynamics}
\dot{x}(t) = f(x(t),u(t)),
\end{equation}
where this system dynamics~(\ref{eq:NominalNonlinearDynamics}) will be used for constructing the equality constraint of the optimization problem. By linearizing the nonlinear system in~(\ref{eq:NominalNonlinearDynamics}) at the equilibrium $(0,0)$, we can obtain the linearized state-space model:
\begin{equation} \label{eq:linearizedDynamics}
\dot{x}(t) = Ax(t)+Bu(t),
\end{equation}
where $A=\frac{\partial f}{\partial x}|_{(0,0)}$ and $B=\frac{\partial f}{\partial u}|_{(0,0)}$.
\par
In the following, we introduce a conventional assumption for the linearized model in~(\ref{eq:linearizedDynamics}), which will be used in the following lemma.
\begin{assumption} \label{ass:a1}
	There exists a state-feedback gain $K$ such that $A+BK$ is stable.
\end{assumption}
\par 
Then, a conventional result regarding the control invariant property of the nonlinear system~(\ref{eq:NominalNonlinearDynamics}) is stated as follows. 
\begin{lemma} \label{lem:chen1998}
	~\cite{chen1998quasi}~If $f:\mathbb{R}^n\times\mathbb{R}^n\to\mathbb{R}^n$ is twice continuously differentiable, $f(0,0)=0$, $u(t)$ is piece-wise right-continuous and Assumption \ref{ass:a1} holds, then given two PD matrices $Q$ and $R$, there exist a state-feedback gain $K$ and a constant $\kappa>0$ such that: 1) The Lyapunov equation $(A+BK+\kappa I)^{\top}P+P(A+BK+\kappa I)=-Q^*$ admits a unique solution $P\succ 0$, where $Q^*=Q+K^{\top}RK\in\mathbb{R}^{n\times n}$ and $\kappa$ is smaller than the real part of $-\overline{\lambda}(A+BK)$; 2) $\Omega({\epsilon}):=\{x\in\mathbb{R}^n|V(x(t))\leq\epsilon\}$ is control invariant by the feedback control law $u(t)=Kx(t)$ for the system in~(\ref{eq:NominalNonlinearDynamics}); 3) $\dot{V}(x(t)) \defeq\normsq{x(t)}_P \leq -\normsq{x(t)}_{Q^*}$ and $u(t)=Kx(t)\in\Ucal$ for $x(t)\in\Omega(\epsilon)$.
\end{lemma}
\par
\begin{figure}
	\centering
	\includegraphics[width=0.8\linewidth]{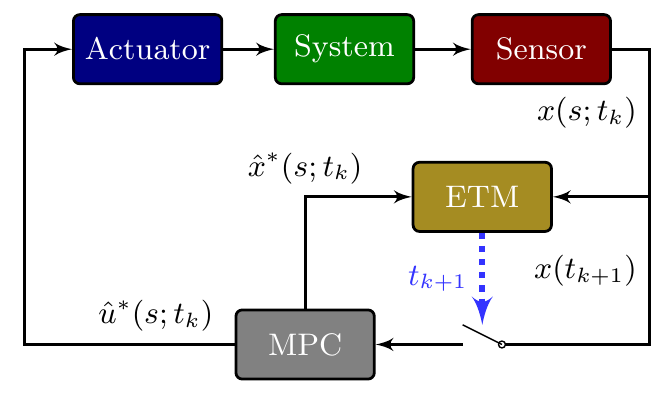}
	\caption{An overview of the event-triggered network control system.}
	\label{fig:diagram_etrhc}
\end{figure}

In the networked control system as illustrated in Fig.~\ref{fig:diagram_etrhc}, the sampling time instants are denoted by $\{t_0,t_1,\ldots,t_k,\ldots\}$ for the nonlinear system in~(\ref{eq:nonlinearDynamics}). At each sampling time instant $t_k$, the MPC controller generates an optimal control sequence $\hat{u}^*(s;t_k)$ and a corresponding optimal predicted state sequence $\hat{x}^*(s;t_k)$ by solving an online optimization problem, where $s\in[t_k,t_k+T]$. In order to save communication resources, we propose to use the event-triggered MPC strategy, where the control inputs are updated aperiodically with larger time intervals rather than periodically with a fixed small time interval. In particular, we aim at designing a more efficient event-triggered MPC scheme for the system in~(\ref{eq:nonlinearDynamics}) such that the better communication performance and the closed-loop stability can be obtained.

\section{Integral-type Event-Triggered MPC} \label{sec:INTET-MPC}

%  The optimal predicted  noticed that these event-triggered instants will be generated by using the integral-type ETM, which will be elaborated after presenting the optimization problem. Then, the resulting closed-loop system is proposed, combining the event-triggered scheme and the MPC , which are generated by using the integral-type ETM

\subsection{Optimization Problem}
In order to avoid ambiguity, we take explicit notations $\hat{u}(s;t_k)$ and $\hat{x}(s;t_k)$ as the control and predicted state sequence at the $\kth$ sampling instant, respectively. Then, the MPC optimization problem at $t_k$ can be designed as
\begin{align} \label{eq:AETDRHC}
\hat{u}^*(s;t_k)=&\arg\min_{\hat{u}\in\mathcal{U}}J\paren{\hat{x}(s;t_k),\hat{u}(s;t_k)} \\
\text{s.t.}\quad& \dot{\hat{x}}(s;t_k) = f(\hat{x}(s;t_k),\hat{u}(s;t_k)) \\
& \hat{u}(s;t_k) \in \mathcal{U}, \quad s \in [t_k,t_k+T] \\
& \norm{\hat{x}(s;t_k)}_{P} \leq \frac{(t_k+T-s)M+s-t_k}{T}\alpha \epsilon \label{eqInAETDRHC:robustconstraint}
\end{align}
where the cost function in (4) is defined by
\begin{equation} \label{eq:costfunction}
\begin{gathered}
J\paren{\hat{x}(s;t_k),\hat{u}(s;t_k)} \defeq \int_{t_k}^{t_k+T}\normsq{\hat{x}(s;t_k)}_{Q}+\normsq{\hat{u}(s;t_k)}_{R} \dopr{s} \\+ \normsq{\hat{x}(t_k+T;t_k)}_{P}.
\end{gathered}
\end{equation}
Note that $Q \succeq 0$, $R \succ 0$, $\MATRIX{I}_n$ is the $n\times n$ identity matrix, $P$ is defined by using the method from Lemma~\ref{lem:chen1998}, $T$ is the prediction horizon, $\epsilon$ is the designed parameter for defining the terminal set, $\alpha\in(0,1)$ is the scaling ratio, and $M$ is the contraction rate for the robustness constraint~(\ref{eqInAETDRHC:robustconstraint}). It is also worthwhile pointing out that the terminal constraint is $\hat{x}(t_k+T;t_k)\in\Omega(\alpha\epsilon)\defeq \{x:\norm{x}_{P} \leq \alpha\epsilon\}$ by the robustness constraint~(\ref{eqInAETDRHC:robustconstraint}). After solving the optimization problem at $t_k$, the optimal control sequence can be expressed as $\hat{u}^*(s;t_k)$, and the corresponding optimal predicted state sequence can be given as $\hat{x}^*(s;t_k)$, where $s\in[t_k,t_k+T]$.
\ifFullVersion
\begin{remark}
	In the optimization problem~(\ref{eq:AETDRHC}), the robustness constraint~(\ref{eqInAETDRHC:robustconstraint}) is used for tightening the $P$-weighted norm of state predictions, which grants MPC capability of compensating the additive disturbance. The robustness constraint-based method for MPC has been firstly proposed in~\cite{li2014event}, where its contraction speed is proportional to the prediction time. In our configuration, the less-conservative robustness constraint is designed with a constant contraction speed as the predicted state trajectory shrink into the terminal set. Compared with the conventional one, our proposed robustness constraint can provide a larger initial feasible set for solving the optimization problem.
\end{remark}
\fi

\subsection{Integral-type Event-Triggered Mechanism}
%For example, when the discrepancy is large enough to reach a specific threshold at some time instant between $t_k$ and $t_k+T]$, then this time instant will be the next sampling instant $t_{k+1}$.
%Note that $\hat{x}^*(s;t_{k})$ is the solution of the differential equation~(\ref{eq:NominalNonlinearDynamics}) by applying the optimal control sequence $\hat{u}^*(s;t_{k})$.
% we take an integral-type ETM for determining when the next optimization should be conducted, i.e. solving the optimization problem and transmitting the optimal control signal to the actuator.
\thirdRevision{In order to more efficiently reduce the frequency of solving the optimization, state sampling, and information transmission from the controller to the actuator, an integral-type ETM is introduced for scheduling and implementing these tasks, i.e. determining the event-triggered instants $\{t_0,t_1,\ldots,t_k,\ldots\}$. At the $\kth$ sampling instant $t_k$, from the system in~(\ref{eq:nonlinearDynamics}) and the MPC optimization problem in~(\ref{eq:AETDRHC}), one can get the real state sequence $x(s;t_k)$ and the optimal predicted state sequence $\hat{x}^*(s;t_k)$ for $s\in[t_k,t_k+T]$, where $T$ denotes the prediction horizon defined in the MPC optimization problem. Due to the additive disturbance in the system, the two sequences $x(s;t_k)$ and $\hat{x}^*(s;t_k)$ cannot coincide with each other. Motivated by this fact, the discrepancy of these two sequences can be used to construct conditions for triggering the next sampling instant.} Therefore, the triggering condition for the integral-type ETM is designed as
\begin{equation} \label{eq:firstTrigger}
\begin{split}
H_{k}&=\inf_{h>0}\left\{h:\int_{t_{k}}^{t_{k}+h}\norm{x(s;t_{k})-\hat{x}^*(s;t_{k})}_{P} \dopr{s} = \delta \right\}, \\
t_{k+1}&=\min\{t_{k}+H_{k},t_{k}+T\},
\end{split}
\end{equation}
\secondRevision{where $h$ is the time variable to be determined for satisfying the equality in the triggering condition~(\ref{eq:firstTrigger}) at $t_k$ and the infimum of $h$ is denoted by $H_k$. If the integral of errors in~(\ref{eq:firstTrigger}) is large enough to reach a specific threshold at some time instant $t_{k}+H_{k}$ between $t_k$ and $t_k+T$, then this $t_{k}+H_{k}$ will be the next sampling instant $t_{k+1}$.}
\begin{remark}
	Note that the designed ETM is based on the integral of errors between actual state sequence and optimal predicted state sequence, which is different from the event-triggered setting in~\cite{li2014event}. The benefit of introducing the integral of errors to the triggering condition lies in that the average of errors between two consecutive event-triggered instants is taken into consideration, leading to some advantages in terms of reducing inter-execution sampling intervals, which will be illustrated in the simulation study in Section~\ref{sec:simulationresult}.
\end{remark}
\par
The following theorem shows some important properties of the proposed integral-type ETM on avoiding the Zeno behavior.
\par 
\begin{theorem} \label{theo:minimumtimeinterval}
	For the nonlinear system in~(\ref{eq:nonlinearDynamics}), if the event-triggered time instants $t_k,\,k\in\mathbb{N}$ are implemented according to~(\ref{eq:firstTrigger}), then the following properties hold: 1) The upper bound on the inter-execution time is $\sup_{k\in\mathbb{N}}(t_{k+1}-t_{k})=T$; 2) the lower bound $\inf_{k\in\mathbb{N}}(t_{k+1}-t_{k})=\beta T$ can be guaranteed by properly designing the triggering level $\delta$ as
	\begin{equation} \label{eqInTh1:minimuminterval}
	\delta = \rho\overline{\lambda}(\sqrt{P})\left[e^{L\beta T}\left(\frac{\beta T}{L}-\frac{1}{L^2}\right)+\frac{1}{L^2}\right],
	\end{equation}
	where $\beta\in(0,1)$ is a scaling parameter.
\end{theorem}
\ifFullVersion
\begin{proof}
	This proof can be done by two steps.
	\par 
	Step 1: The upper bound of inter-execution intervals is $T$. From the design of the integral-type ETM, it can be directly deduced that all the intervals $t_{k+1}-t_k$ is less than or equal to the prediction horizon $T$.
	\par 
	Step 2: The inter-execution intervals can be lower bounded to $\beta T$ by properly designing the triggering level $\delta$ according to~(\ref{eqInTh1:minimuminterval}). To prove this result, we firstly consider the upper bound for $\norm{x(s;t_k)-\hat{x}^*(s;t_k)}_{P}$ at $t_k$. We assume here that the sensor measurements are accurate. Thus, it follows that $x(t_k;t_k)-\hat{x}^*(t_k;t_k)=0$. By using the triangle inequality, we have
	$
	\norm{x(s;t_k)-\hat{x}^*(s;t_k)}_{P} 
	\leq\norm{x(t_k;t_k)-\hat{x}^*(t_k;t_k)+\int_{t_k}^{s}\dot{x}(\tau;t_k)-\dot{x}^*(\tau;t_k)\dopr{\tau}}_{P} 
	\leq\int_{t_k}^{s}L\norm{x(\tau;t_k)-\hat{x}^*(\tau;t_k)}_{P}+\norm{\omega(\tau)}_{P}\dopr{\tau} 
	\leq\int_{t_k}^{s}L\norm{x(\tau;t_k)-\hat{x}^*(\tau;t_k)}_{P}\dopr{\tau}+\rho\overline{\lambda}(\sqrt{P})(s-t_k).	   
	$
	Then by applying the integral form of Gronwall-Bellman inequality, it can be obtained that 
	$
	\norm{x(s;t_k)-\hat{x}^*(s;t_k)}_{P} \leq \rho\overline{\lambda}(\sqrt{P})(s-t_k)e^{L(s-t_k)}.
	$
	Substituting the previous inequality to~(\ref{eq:firstTrigger}), we can deduce that
	$
	\int_{t_k}^{t_{k+1}}\norm{x(s;t_k)-\hat{x}^*(s;t_k)}_{P} \dopr{s}  \leq\int_{t_k}^{t_{k+1}}\rho\overline{\lambda}(\sqrt{P})(s-t_k)e^{L(s-t_k)}\dopr{s} 
	=\rho \overline{\lambda} (\sqrt{P})\left[e^{L(t_{k+1}-t_k)}\left(\frac{t_{k+1}-t_k}{L}-\frac{1}{L^2}\right)+\frac{1}{L^2}\right].
	$ Since $\overline{\lambda}(\sqrt{P})(s-t_k)e^{L(s-t_k)}$ is strictly larger than zero for $s>t_k$, we can choose $\delta$ as~(\ref{eqInTh1:minimuminterval}) such that the lower bound of the triggered time interval is $\inf_{k\in\mathbb{N}}(t_{k+1}-t_{k})=\beta T$. The proof is completed.
		%	\begin{equation*}
	%	\begin{split} 
	%	&\int_{t_k}^{t_{k+1}}\norm{x(s;t_k)-\hat{x}^*(s;t_k)}_{P} \dopr{s} \\ \leq&\int_{t_k}^{t_{k+1}}\frac{\rho\overline{\lambda}(\sqrt{P})(s-t_k)e^{L(s-t_k)}}{t_{k+1}-t_k}\dopr{s} \\
	%	=&\rho \overline{\lambda} (\sqrt{P})\left[e^{L(t_{k+1}-t_k)}\left(\frac{1}{L}-\frac{1}{(t_{k+1}-t_k)L^2}\right)+\frac{1}{(t_{k+1}-t_k)L^2}\right].
	%	\end{split}   
	%	\end{equation*}
\end{proof}
\else
\begin{proof}
	The proof can be completed by following the similar methodology as for Theorem 3 in~\cite{li2014event}.
\end{proof}
\fi

\subsection{The Closed-Loop System}

By implementing the generated optimal control sequence into the system model in~(\ref{eq:nonlinearDynamics}), we write the resulting closed-loop system as
\begin{equation} \label{eq:closedloopSystem}
\dot{x}(t) = f(x(t),\hat{u}^*(t;t_{k})) + \omega(t), \quad k\in\{0,1,2,\ldots\},
\end{equation}
where $t_{k}$ is the $\kth$ event-triggered instant generated by~(\ref{eq:firstTrigger}). The feasibility of MPC problem~(\ref{eq:AETDRHC}) and the closed-loop stability of~(\ref{eq:closedloopSystem}) will be respectively analyzed in the following section. For a clear view of the aforementioned integral-type ET-MPC, we design the ET-MPC algorithm as described in Algorithm~\ref{alg:int_ETRHC}.
\begin{algorithm}
	\caption{Integral-type ET-MPC}
	\label{alg:int_ETRHC}
	\begin{algorithmic}[1]
		\WHILE{The control action is not stopped}
		\IF{$k$ = $0$}
		\STATE Solve the optimization problem in~(\ref{eq:AETDRHC}) at $t_0$;
		\ENDIF
		\WHILE {The ETM condition~(\ref{eq:firstTrigger}) is not triggered}
		\STATE Apply optimal control input $\hat{u}^*(s;t_k)$, where $s\in[t_k,t_k+T]$;
		\ENDWHILE
		\STATE Obtain $t_{k+1}$ according to~(\ref{eq:firstTrigger});
		\STATE Solve the optimization problem in~(\ref{eq:AETDRHC}) at $t_{k+1}$;
		\ENDWHILE
	\end{algorithmic}
\end{algorithm}

\section{Main Results} \label{sec:mainresult}
\subsection{Feasibility Analysis}
Following a conventional setup for MPC framework, we construct a classical feasible control sequence for the optimization problem~(\ref{eq:AETDRHC}) as follows:
\begin{equation} \label{eq:dualmodeControl}
\tilde{u}(s;t_k)=\left\{
\begin{aligned}
&\hat{u}^*(s;t_{k-1}), &\text{if } s \in [t_{k},t_{k-1}+T] \\
&K x(s;t_k), &\text{if } \thirdRevision{s \in (t_{k-1}+T,t_{k}+T]}\\
\end{aligned}
\right.
\end{equation}
Then the candidate state sequence evolves as
\begin{equation} \label{eq:closedLoop}
\dot{\tilde{x}}(s;t_k)=f(\tilde{x}(s;t_k),\tilde{u}(s;t_k)).
\end{equation}
Before presenting the result of this section, we introduce a lemma that will be useful in the following analysis.
\begin{lemma} \label{lem:lem1}
	\secondRevision{Let $g: \mathbb{R} \to \mathbb{R}^N$ be a differentiable vector-valued function defined on $t\in[a,b]$}. Then the following inequality holds:
	\begin{equation} \label{eq:lemresult}
	\sup_{t\in[a,b]}\norm{g(t)} \leq \frac{1}{2}\int_{a}^{b}\norm{g^{'}(t)}\dopr{t}+\frac{1}{2}\norm{g(a)+g(b)},
	\end{equation}
	\secondRevision{where $g^{'}$ is the derivative of $g$}.
\end{lemma}
\ifFullVersion
\begin{proof}
	For every $t\in[a,b]$, we have two results:
$
	g(t)=g(a)+\int_{a}^{t}g^{'}(\tau)\dopr{\tau}, \label{eq:lem11}\quad
	g(b)=g(t)+\int_{t}^{b}g^{'}(\tau)\dopr{\tau}. 
$
	Subtracting the aforementioned two equations yields
$
	2g(t)=g(a)+g(b)+\int_{a}^{t}g^{'}(\tau)\dopr{\tau}+\int_{b}^{t}g^{'}(\tau)\dopr{\tau}.
$
	By employing the triangle inequality, we can deduce from the above equality that 
$
	\norm{g(t)}\leq \frac{1}{2}\norm{g(a)+g(b)}+\frac{1}{2}\int_{a}^{t}\norm{g^{'}(\tau)}\dopr{\tau}+\frac{1}{2}\int_{t}^{b}\norm{g^{'}(\tau)}\dopr{\tau}
	=\frac{1}{2}\int_{a}^{b}\norm{g^{'}(t)}\dopr{t}+\frac{1}{2}\norm{g(a)+g(b)}.
$
	Since every $g(t)$ defined on the closed interval $[a,b]$ is equal to or less than the right side of the above inequality, the result in~(\ref{eq:lemresult}) thus holds.
\end{proof}
\else
\begin{proof}
	Due to the space limit, please refer to the detailed technical proof in our full version~\cite{sun2018integral}.
\end{proof}
\fi
\par 
To facilitate the analysis of integral-type event-triggered configuration, we make use of a term $e(s;t_k)=\tilde{x}(s;t_k)-\hat{x}^{*}(s;t_{k-1})$ for $s\in[t_k,t_{k-1}+T]$. Note that $\tilde{x}(s;t_k)$ is defined as the candidate state sequence generated by the nominal system~(\ref{eq:closedLoop}), and $\hat{x}^{*}(s;t_{k-1})$ is the solution of the optimization problem~(\ref{eq:AETDRHC}). In addition, we let $\tilde{x}(t_k;t_k)=x(t_k)$ by sampling the state at the $\kth$ triggered time instant $t_k$, since we assume that there is no measurement inaccuracy.
%Note from Lipschitz condition that we have
%\begin{gather*}
%\dot{e}(s;t_k)=f(\tilde{x}(s;t_k),\tilde{u}(s;t_k))-f(\hat{x}(s;t_k),\hat{u}(s;t_k))  \\
%e(s;t_k)\leq \rho(s-t_k) e^{L(s-t_k)} \\
%e^{\top}(s;t_k)Pe(s;t_k)=\norm{\tilde{x}(s;t_k)-\hat{x}(s;t_k)}_{P}^2 \\
%\norm{\tilde{x}(s;t_k)-\hat{x}(s;t_k)}_{P} = \norm{e(s;t_k)}_{P} \leq \rho \bar{\lambda}(\sqrt{P})(s-t_{k})e^{L(s-t_{k})}
%\end{gather*}
Then we can propose the following result.
\begin{corollary} \label{cor:cor1}
	Given $e(s;t_k)$ defined on $s\in[t_k,t_{k-1}+T]$, the following inequality holds
	\begin{equation} \label{eq:cor1result}
		\sup_{s\in[t_k,t_{k-1}+T]}\norm{e(s;t_k)}_{P}\leq \frac{L^2\beta T}{L\beta T-1} e^{L(1-\beta)T}\delta,
		%\leq &(L+\frac{L^2\beta T}{L\beta T-1})\int_{t_k}^{t_{k-1}+T} \norm{e(s;t_k)}_{P}\dopr{s}.
	\end{equation}
	under the condition that $L\beta T > 1$.
\end{corollary}
\ifFullVersion
\begin{proof}
	By Lemma~\ref{lem:lem1}, we can obtain that 	
	%$
	%		\sup_{s\in[t_k,t_{k-1}+T]}\norm{\sqrt{P}e(s;t_k)} 
	%		\leq \frac{1}{2}L\int_{t_k}^{t_{k-1}+T} \norm{\sqrt{P}e(s;t_k)}\dopr{s} + \frac{1}{2}\norm{\sqrt{P}e(t_k;t_k)+\sqrt{P}e(t_{k-1}+T;t_k)}.
	%$
	\begin{equation}
	\begin{split}
		\sup_{s\in[t_k,t_{k-1}+T]}&\norm{\sqrt{P}e(s;t_k)} 
		\leq \frac{1}{2}L\int_{t_k}^{t_{k-1}+T} \norm{\sqrt{P}e(s;t_k)}\dopr{s}\\ &+ \frac{1}{2}\norm{\sqrt{P}e(t_k;t_k)+\sqrt{P}e(t_{k-1}+T;t_k)}.
	\end{split}
	\end{equation}
		Note that $\norm{e(s;t_k)}_{P}=\norm{\sqrt{P}e(s;t_k)}$. Then it can be deduced from the above inequality that 
	%$
	%	\sup_{s\in[t_k,t_{k-1}+T]}\norm{e(s;t_k)}_{P} 
	%	\leq \frac{1}{2}L\int_{t_k}^{t_{k-1}+T} \norm{e(s;t_k)}_{P}\dopr{s} + \frac{1}{2}\norm{2e(t_{k};t_k) +\int_{t_k}^{t_{k-1}+T}\dot{e}(s;t_k)\dopr{s}}_{P} 
	%	\leq L\int_{t_k}^{t_{k-1}+T} \norm{e(s;t_k)}_{P}\dopr{s} + \norm{e(t_{k};t_k)}_{P}.
	%$
	\begin{equation}\label{eqInCor1:2}
	\begin{split}
	\sup_{s\in[t_k,t_{k-1}+T]}&\norm{e(s;t_k)}_{P} 
	\leq \frac{1}{2}L\int_{t_k}^{t_{k-1}+T} \norm{e(s;t_k)}_{P}\dopr{s}\\ &+ \frac{1}{2}\norm{2e(t_{k};t_k) +\int_{t_k}^{t_{k-1}+T}\dot{e}(s;t_k)\dopr{s}}_{P} 
	\\&\leq L\int_{t_k}^{t_{k-1}+T} \norm{e(s;t_k)}_{P}\dopr{s} + \norm{e(t_{k};t_k)}_{P}.
	\end{split}
	\end{equation}
		Next we show that the upper bound of $\norm{e(t_{k};t_k)}_{P}$ is related to the triggering level $\delta$. By using the result from Theorem~\ref{theo:minimumtimeinterval}, the upper bound of $\norm{x(s;t_{k-1})-\hat{x}^*(s;t_{k-1})}_{P}$ is $\rho\overline{\lambda}(\sqrt{P})(s-t_{k-1})e^{L(s-t_{k-1})}$ for $s\in[t_{k-1},t_{k}]$. Suppose that there exists a maximum value of $\norm{x(s;t_{k-1})-\hat{x}^*(s;t_{k-1})}_{P}$ such that $\norm{x(t_{k-1}+h;t_{k-1})-\hat{x}^*(t_{k-1}+h;t_{k-1})}_{P}\geq\norm{x(s;t_{k-1})-\hat{x}^*(s;t_{k-1})}_{P}$ for $s\in[t_{k-1},t_{k}]$. Here $h\in[\beta T,T]$ is some positive real number. Note from Theorem~\ref{theo:minimumtimeinterval} that we also have $\delta=\rho \overline{\lambda} (\sqrt{P})\left[e^{Lh}\left(\frac{h}{L}-\frac{1}{L^2}\right)+\frac{1}{L^2}\right]$, where the triggering level $\delta$ is designed as the integral of $\norm{x(s;t_{k-1})-\hat{x}^*(s;t_{k-1})}_{P}$ from $t_{k-1}$ to $t_{k-1}+h$. By following some calculation, we can obtain that
	%$
	%	\rho \overline{\lambda} (\sqrt{P})\left[e^{Lh}(h-\frac{1}{L})+\frac{1}{L}\right] = L\delta
	%$, 
	\begin{equation}
	\rho \overline{\lambda} (\sqrt{P})\left[e^{Lh}(h-\frac{1}{L})+\frac{1}{L}\right] = L\delta,
	\end{equation}
	and consequently it follows that
	%$
	% 	\norm{e(h;t_k)}_{P}\leq\rho \overline{\lambda} (\sqrt{P})e^{Lh}h \leq L\delta\frac{h}{h-\frac{1}{L}} 
	%$.
	\begin{equation}
		\sigma(t_{k-1}+h;t_{k-1}) \leq \rho \overline{\lambda} (\sqrt{P})e^{Lh}h \leq L\delta\frac{h}{h-\frac{1}{L}}, 
	\end{equation}
	where $\sigma(t_{k-1}+h;t_{k-1}) \defeq \norm{x(t_{k-1}+h;t_{k-1})-\hat{x}^*(t_{k-1}+h;t_{k-1})}_{P}$ and $L\beta T>1$. Since the function $\frac{h}{h-\frac{1}{L}}$ gets its maximum at $h=\beta T$, the above inequality becomes
	%$
	%		\norm{e(h;t_k)}_{P} \leq \frac{L^2\beta T}{L\beta T-1}\delta
	%$.
	\begin{equation}\label{eqInCor1:last}
	\norm{e(t_k;t_k)}_{P} \leq \sigma(t_{k-1}+h;t_{k-1}) \leq \frac{L^2\beta T}{L\beta T-1}\delta.
	\end{equation}
	According to Gronwall-Bellman inequality, one can obtain~(\ref{eq:cor1result}) by substituting~(\ref{eqInCor1:last}) to~(\ref{eqInCor1:2}). The proof is thus completed.
\end{proof}
\else
\begin{proof}
	Due to the space limit, please refer to the detailed technical proof in our full version~\cite{sun2018integral}.
\end{proof}
\fi
\par 
Now we can analyze the iterative feasibility of the MPC problem~(\ref{eq:AETDRHC}), implying that if the MPC problem admits a solution at current time instant then a feasible solution exists for the next time instant. To prove this result, we use a conventional feasible control sequence candidate $\tilde{u}(s;t_k)$ at time instant $t_k$ defined in~(\ref{eq:dualmodeControl}), where $\tilde{u}(s;t_k)=\hat{u}^{*}(s;t_{k-1})$ for $s\in[t_{k},t_{k-1}+T]$ and $\tilde{u}(s;t_k)=K\tilde{x}(s;t_k)$ for $s\in(t_{k-1}+T,t_{k}+T]$. In the following theorem, we will show that the designed control sequence candidate $\tilde{u}(s;t_k)$ can steer the feasible state $\tilde{x}(s;t_k)$ into $\Omega{(\alpha\epsilon)}$ if some conditions can be satisfied. In addition, it is also necessary to show that the candidate state $\tilde{x}(s;t_{k})$ will fulfill the designed state constraint~(\ref{eqInAETDRHC:robustconstraint}). 
\begin{assumption} \label{ass:a2}
	The optimization problem~(\ref{eq:AETDRHC}) admits a feasible solution $\hat{u}^*(s;t_0)$ at the initial time $t_0$.
\end{assumption}
\begin{theorem}[Feasibility] \label{theo:feasibilityResult}
	Suppose that Assumptions~\ref{ass:a1} and~\ref{ass:a2} hold. The MPC problem~(\ref{eq:AETDRHC}) is iteratively feasible under the following conditions:
	\begin{gather}
		\resizebox{\linewidth}{!}{$
		\frac{L^2\beta Te^{L(1-\beta)T}}{L\beta T-1}\rho\overline{\lambda}(\sqrt{P})\left[e^{L\beta T}\left(\frac{\beta T}{L}-\frac{1}{L^2}\right)+\frac{1}{L^2}\right] \leq (1-\alpha)\epsilon$}, \label{eqInTheo2:c1} \\
		T\geq \max\left\{ -2\frac{\overline{\lambda}(P)}{\underline{\lambda}(Q^*)\beta}\ln{\alpha}, \frac{1}{L\beta} \right\}, \label{eqInTheo2:c2} \\
		M\geq\max\left\{\frac{L^2\beta Te^{L(1-\beta)T}}{L\beta T-1}\frac{\delta}{\alpha\epsilon}+1,1-\frac{1}{\beta}+\frac{1}{\alpha\beta}\right\}. \label{eqInTheo2:c3}
	\end{gather}
	Moreover, the maximum allowable disturbance can be given as
	\begin{equation} \label{eqInTheo2:maxdisturbance}
		\rho\leq \frac{(1-\alpha)\epsilon}{\frac{L^2\beta Te^{L(1-\beta)T}}{L\beta T-1}\overline{\lambda}(\sqrt{P})\left[e^{L\beta T}\left(\frac{\beta T}{L}-\frac{1}{L^2}\right)+\frac{1}{L^2}\right]}. 	   
	\end{equation}
\end{theorem}
\ifFullVersion
\begin{proof}
	First, we show that the designed control sequence $\tilde{u}(s;t_{k})$ for $s\in[t_{k},t_{k-1}+T]$ drives $\tilde{x}(s;t_k)$ into $\Omega{(\epsilon)}$, i.e. $\norm{\tilde{x}(s;t_k)}_{P}\leq\epsilon$. Let us construct an error norm $\norm{\tilde{x}(s;t_k)-\hat{x}^{*}(s;t_{k-1})}_{P}$ for $s\in[t_k,t_{k-1}+T]$. By using Corollary~\ref{cor:cor1}, we can obtain that %$\sup_{s\in[t_k,t_{k-1}+T]}\norm{\tilde{x}(s;t_k)-\hat{x}^{*}(s;t_{k-1})}_{P} \leq \frac{L^2\beta T}{L\beta T-1} e^{L(1-\beta)T}\delta$.
	\begin{equation}
	\begin{split}
		\sup_{s\in[t_k,t_{k-1}+T]}&\norm{\tilde{x}(s;t_k)-\hat{x}^{*}(s;t_{k-1})}_{P} \\ \leq & \frac{L^2\beta T}{L\beta T-1} e^{L(1-\beta)T}\delta,
	\end{split}
	\end{equation}
	which turns out that $\norm{\tilde{x}(t_{k-1}+T;t_k)-\hat{x}^{*}(t_{k-1}+T;t_{k-1})}_{P} \leq \frac{L^2\beta Te^{L(1-\beta)T}}{L\beta T-1} \delta$. With help of the Triangle inequality, we have
	%$\norm{\tilde{x}(t_{k-1}+T;t_k)}_{P}\leq\norm{\hat{x}^{*}(t_{k-1}+T;t_{k-1})}_{P}+\frac{L^2\beta Te^{L(1-\beta)T}}{L\beta T-1}\delta$, 
	\begin{equation}
	\begin{split}
		&\norm{\tilde{x}(t_{k-1}+T;t_k)}_{P} \\ \leq&\norm{\hat{x}^{*}(t_{k-1}+T;t_{k-1})}_{P}+\frac{L^2\beta Te^{L(1-\beta)T}}{L\beta T-1}\delta,
	\end{split}
	\end{equation}
	which implies that $\norm{\tilde{x}(t_{k-1}+T;t_k)}_{P}\leq\alpha\epsilon+\frac{L^2\beta Te^{L(1-\beta)T}}{L\beta T-1}\delta$.
	\par 
    Note from Theorem~\ref{theo:minimumtimeinterval} that the triggering level is designed as $\delta = \rho\overline{\lambda}(\sqrt{P})\left[e^{L\beta T}\left(\frac{\beta T}{L}-\frac{1}{L^2}\right)+\frac{1}{L^2}\right]$. In order to steer the candidate state sequence $\tilde{x}(t_{k-1}+T;t_k)$ into $\Omega{(\epsilon)}$, one can simply deduce that the following inequality must hold:
	\begin{equation} \label{eqInTheo2:disturbanceBound}
	\begin{split}
		&\frac{L^2\beta Te^{L(1-\beta)T}}{L\beta T-1}\rho\overline{\lambda}(\sqrt{P})\left[e^{L\beta T}\left(\frac{\beta T}{L}-\frac{1}{L^2}\right)+\frac{1}{L^2}\right] \\ \leq & (1-\alpha)\epsilon.
	\end{split}
	\end{equation} 
    From~(\ref{eqInTheo2:disturbanceBound}), it can be also obtained that the maximum bound for disturbance satisfies~(\ref{eqInTheo2:maxdisturbance}).
    \par 
	Second, we consider the candidate sequence $\tilde{x}(s;t_k)$ for $s\in(t_{k-1}+T,t_k+T]$. Then using Lemma~\ref{lem:chen1998}, we can verify that $\Omega(\epsilon)$ is an invariant set for the closed-loop system $\dot{\tilde{x}}(s;t_k)=f\left(\tilde{x}(s;t_k),K(\tilde{x}(s;t_k))\right)$. Consequently, we can deduce that $\dot{V}(\tilde{x}(s;t_k))\leq -\normsq{\tilde{x}(s;t_k)}_{Q^*}$. By the virtue of comparison principle, it follows that
	\begin{equation} \label{eqInTheo2:feedbackX}
		V(\tilde{x}(s;t_k))\leq \epsilon^2e^{-\frac{\underline{\lambda}(Q^*)}{\overline{\lambda}(P)}(s-t_{k-1}-T)}
	\end{equation}
	for $s\in(t_{k-1}+T,t_k+T]$, which indicates that $V(\tilde{x}(t_{k}+T;t_k))\leq \epsilon^2e^{-\frac{\underline{\lambda}(Q^*)}{\overline{\lambda}(P)}(t_k-t_{k-1})}$. By using Theorem~\ref{theo:minimumtimeinterval}, we can have $\inf(t_k-t_{k-1})=\beta T$. To obtain $\norm{\tilde{x}(t_{k}+T;t_k)}_P\leq\alpha\epsilon$, it is equivalent to show that $V(\tilde{x}(t_{k}+T;t_k))\leq\alpha^2\epsilon^2$. With some calculation, one can have $T\geq -2\frac{\overline{\lambda}(P)}{\underline{\lambda}(Q^*)\beta}\ln{\alpha}$ such that the previous inequality holds. Similar argument can be found in~\cite{li2014event}.
	\par 
	Third, we show that $\tilde{x}(s;t_{k})$ satisfies the state constraint~(\ref{eqInAETDRHC:robustconstraint}). For $s\in(t_{k},t_{k-1}+T]$, one can get %$\norm{\tilde{x}(s;t_{k})}_P\leq \norm{\hat{x}^{*}(s;t_{k-1})}_P+\frac{L^2\beta Te^{L(1-\beta)T}}{L\beta T-1}\delta$
	\begin{equation}
		\norm{\tilde{x}(s;t_{k})}_P\leq \norm{\hat{x}^{*}(s;t_{k-1})}_P+\frac{L^2\beta Te^{L(1-\beta)T}}{L\beta T-1}\delta,
	\end{equation}
	which can be derived from~(\ref{eq:cor1result}). Then we need to prove %$\frac{(t_{k}+T-s)M+s-t_{k}}{T}\alpha\epsilon \leq \frac{(t_{k-1}+T-s)M+s-t_{k-1}}{T}\alpha \epsilon + \frac{L^2\beta Te^{L(1-\beta)T}}{L\beta T-1}\delta$.
	\begin{equation}
	\begin{split}
		&\frac{(t_{k}+T-s)M+s-t_{k}}{T}\alpha\epsilon \\ \leq & \frac{(t_{k-1}+T-s)M+s-t_{k-1}}{T}\alpha \epsilon + \frac{L^2\beta Te^{L(1-\beta)T}}{L\beta T-1}\delta.
	\end{split}
	\end{equation} 
	By some calculation, it can be obtain that $M\geq \frac{L^2\beta Te^{L(1-\beta)T}}{L\beta T-1}\frac{\delta}{\alpha\epsilon}+1$. For $s\in(t_{k-1}+T,t_{k}+T]$, it can be deduced from~(\ref{eqInTheo2:feedbackX}) that 
	%$\norm{\tilde{x}(s;t_{k})}_P \leq \epsilon e^{-\frac{\underline{\lambda}(Q^*)}{\overline{\lambda}(P)}(s-t_{k-1}-T)/2}$. 
	\begin{equation}
		\norm{\tilde{x}(s;t_{k})}_P \leq \epsilon e^{-\frac{\underline{\lambda}(Q^*)}{\overline{\lambda}(P)}(s-t_{k-1}-T)/2}.
	\end{equation}
	In order to prove $\norm{\tilde{x}(s;t_{k})}_P\leq \frac{(t_{k}+T-s)M+s-t_{k}}{T}\alpha\epsilon$, it is equivalent to show 
	%$\frac{(t_{k}+T-s)M+s-t_{k}}{T}\alpha\epsilon \geq \epsilon e^{-\frac{\underline{\lambda}(Q^*)}{\overline{\lambda}(P)}(s-t_{k-1}-T)/2}$.
	\begin{equation}
		\frac{(t_{k}+T-s)M+s-t_{k}}{T}\alpha\epsilon \geq \epsilon e^{-\frac{\underline{\lambda}(Q^*)}{\overline{\lambda}(P)}(s-t_{k-1}-T)/2}.
	\end{equation}
	For brevity, we denote $F(s)=\frac{T/\alpha \cdot e^{-\frac{\underline{\lambda}(Q^*)}{\overline{\lambda}(P)}(s-t_{k-1}-T)/2}+t_{k}-s}{t_{k}+T-s}$, and it turns out that $M\geq F(s)$.  By evaluating the derivative of $F(s)$, it can be verified that $F^{'}(s)$ is non-positive for $s\in(t_{k-1}+T,t_{k}+T]$, which indicates $M\geq 1-\frac{1}{\beta}+\frac{1}{\alpha\beta}$. Finally, the designing parameter should be configured as $M\geq\max\{\frac{L^2\beta Te^{L(1-\beta)T}}{L\beta T-1}\frac{\delta}{\alpha\epsilon}+1,1-\frac{1}{\beta}+\frac{1}{\alpha\beta}\}$ for guaranteeing the satisfaction of the proposed robustness constraint.
	The proof is completed.
\end{proof}
\else
\begin{proof}
	Due to the space limit, please refer to the detailed technical proof in our full version~\cite{sun2018integral}.
\end{proof}
\fi
\begin{remark}
	Note from Theorem~\ref{theo:feasibilityResult} that the feasibility can be affected by the prediction horizon $T$, the Lipschitz constant $L$, and the disturbance bound $\rho$. In order to achieve the recursive feasibility, the prediction horizon $T$ and the design parameter $M$ in~(\ref{eqInAETDRHC:robustconstraint}) should be both lower bounded. It should be also noted that the maximum allowable disturbance bound can be estimated according to~(\ref{eqInTheo2:maxdisturbance}).
\end{remark}
\subsection{Stability Analysis}
In this part, we investigate the closed-loop stability by applying the proposed integral-type ET-MPC. It is worthwhile to point out that, due to the disturbance, the established  closed-loop stability property can steer system states into an invariant set. Under the MPC configuration, the analysis for stability can be divided into two steps: 1) The first step is to ensure that the system trajectory will enter the terminal set in finite time; 2) the second step is to prove that the closed-loop system is stable after the state enters the terminal set.
\begin{theorem}[Stability] \label{theo:stability}
	Suppose that Assumptions~\ref{ass:a1} and~\ref{ass:a2} hold, and the conditions in Theorem~\ref{theo:feasibilityResult} are satisfied. The state of the closed-loop system~(\ref{eq:closedloopSystem}) enters the designed terminal set in finite time and converges to $\Omega{(\bar{\epsilon})}$ if the following condition holds for some $n\in\mathbb{N}$:
	\begin{equation} \label{cond:stability1}
		\begin{split}
			&\resizebox{\linewidth}{!}{$\frac{\overline{\lambda}(Q)}{\underline{\lambda}(P)}\frac{L^2(1-\beta)T}{L\beta T-1}e^{L(1-\beta)T}\delta \bigg[ \frac{L^2\beta T}{L\beta T-1}e^{L(1-\beta)T}\delta
			 + 2[(1-\beta)M+\beta]\alpha\epsilon \bigg]$} \\& \resizebox{\linewidth}{!}{$ + \frac{L^4\beta T}{(L\beta T-1)^2} e^{2L(1-\beta)T}\delta^2 \leq \frac{\underline{\lambda}(Q)n}{\overline{\lambda}(P)(n+1)} \left(\alpha\epsilon-\frac{L^2\beta T}{L\beta T-1}\delta\right)^2$}.
		\end{split}
	\end{equation}
\end{theorem}
\ifFullVersion
\begin{proof}
	This theorem will be proved by two steps. 
	\par
	Step 1: For all initial state $x(t_0)\in \mathcal{X}\setminus\Omega{(\alpha\epsilon)}$ where $\mathcal{X}$ is the initial feasible set for system state, we aim to show the state trajectory enters $\Omega{(\alpha\epsilon)}$ in finite time. In this situation, we construct an error term as $\Delta \tilde{J}\paren{x(s;t_k),u(s;t_k)}:=J\paren{\tilde{x}(s;t_k),\tilde{u}(s;t_k)}  -J\paren{\hat{x}^{*}(s;t_{k-1}),\hat{u}^{*}(s;t_{k-1})}$. Expanding this term yields 
	\begin{equation}
	\begin{split}
		&\Delta \tilde{J}\paren{x(s;t_k),u(s;t_k)} \\ &= \int_{t_{k}}^{t_{k}+T} \normsq{\tilde{x}(s;t_k)}_Q + \normsq{\tilde{u}(s;t_k)}_R\dopr{s}+\normsq{\tilde{x}(t_k+T;t_k)}_P \\& \qquad - \int_{t_{k-1}}^{t_{k-1}+T} \normsq{\hat{x}^*(s;t_{k-1})}_Q + \normsq{\hat{u}^*(s;t_{k-1})}_R\dopr{s} \\ & \qquad -\normsq{\hat{x}^*(t_{k-1}+T;t_{k-1})}_P.
	\end{split}
	\end{equation}
	Substituting $\tilde{u}(s;t_k)$ in~(\ref{eq:dualmodeControl}) to the above equation, we can obtain that
	\begin{equation}
	\begin{split}
		\Delta \tilde{J}&\paren{x(s;t_k),u(s;t_k)} = \int_{t_{k-1}+T}^{t_{k}+T} \normsq{\tilde{x}(s;t_k)}_{Q^*} \dopr{s}\\&+\int_{t_{k}}^{t_{k-1}+T} \normsq{\tilde{x}(s;t_k)}_Q-\normsq{\hat{x}^*(s;t_{k-1})}_Q \\&\qquad\qquad+ \normsq{\tilde{u}(s;t_k)}_R-\normsq{\hat{u}^*(s;t_{k-1})}_R\dopr{s}\\&-\int_{t_{k-1}}^{t_k}\normsq{\hat{x}^*(s;t_{k-1})}_Q + \normsq{\hat{u}^*(s;t_{k-1})}_R\dopr{s} \\ &+ \normsq{\tilde{x}(t_k+T;t_k)}_P - \normsq{\hat{x}^*(t_{k-1}+T;t_{k-1})}_P.
	\end{split}
	\end{equation}
	Note from Lemma~\ref{lem:chen1998} that $\dot{V}(\tilde{x}(s;t_k))\leq -\normsq{\tilde{x}(s;t_k)}_{Q^*}$. Taking integral from $t_{k-1}+T$ to $t_{k}+T$ of the above inequality yields 
	\begin{equation}
	\begin{split}
		&\int_{t_{k-1}+T}^{t_{k}+T}\dot{V}(\tilde{x}(s;t_k))\dopr{s} \\
		= &\normsq{\tilde{x}(t_k+T;t_k)}_P - \normsq{\tilde{x}(t_{k-1}+T;t_{k})}_P \\ 
		\leq& -\int_{t_{k-1}+T}^{t_{k}+T}\normsq{\tilde{x}(s;t_k)}_{Q^*}\dopr{s}.
	\end{split}
	\end{equation}
	Applying this fact to $\Delta \tilde{J}\paren{x(s;t_k),u(s;t_k)}$, it can be shown that 
	\begin{equation} \label{eqinTheo:3}
	\begin{split} 
	&\Delta \tilde{J}\paren{x(s;t_k),u(s;t_k)} \\\leq& \underbrace{\int_{t_{k}}^{t_{k-1}+T} \normsq{\tilde{x}(s;t_k)}_Q-\normsq{\hat{x}^*(s;t_{k-1})}_Q \dopr{s}}_{\Delta \tilde{J}_1} \\& \underbrace{-\int_{t_{k-1}}^{t_k}\normsq{\hat{x}^*(s;t_{k-1})}_Q + \normsq{\hat{u}^*(s;t_{k-1})}_R\dopr{s}}_{\Delta \tilde{J}_2} \\ & + \underbrace{\normsq{\tilde{x}(t_{k-1}+T;t_{k})-\hat{x}^*(t_{k-1}+T;t_{k-1})}_P}_{\Delta \tilde{J}_3}
	\end{split}
	\end{equation} 
	To analyze the above inequality, we firstly consider the term 
	%$A=\int_{t_{k}}^{t_{k-1}+T} \normsq{\tilde{x}(s;t_k)}_Q-\normsq{\hat{x}^*(s;t_{k-1})}_Q \dopr{s}$.
	\begin{equation}
		\Delta \tilde{J}_1=\int_{t_{k}}^{t_{k-1}+T} \normsq{\tilde{x}(s;t_k)}_Q-\normsq{\hat{x}^*(s;t_{k-1})}_Q \dopr{s}.
	\end{equation}
	By using the triangle inequality, we have 
	%$A \leq \int_{t_{k}}^{t_{k-1}+T} \normsq{\tilde{x}(s;t_k)-\hat{x}^*(s;t_{k-1})}_Q \dopr{s}+2\int_{t_{k}}^{t_{k-1}+T} \norm{\tilde{x}(s;t_k)-\hat{x}^*(s;t_{k-1})}_Q\cdot\norm{\hat{x}^*(s;t_{k-1})}_Q \dopr{s}$.
	\begin{equation}
	\begin{split}
		\Delta \tilde{J}_1 \leq &\int_{t_{k}}^{t_{k-1}+T} \normsq{\tilde{x}(s;t_k)-\hat{x}^*(s;t_{k-1})}_Q \dopr{s} \\ &+2\int_{t_{k}}^{t_{k-1}+T} \norm{\tilde{x}(s;t_k)-\hat{x}^*(s;t_{k-1})}_Q \\ &\qquad\quad\qquad\cdot\norm{\hat{x}^*(s;t_{k-1})}_Q \dopr{s}
	\end{split}
	\end{equation}
	Then apply Holder inequality, and it follows that 
	%$A \leq \int_{t_{k}}^{t_{k-1}+T} \norm{\tilde{x}(s;t_k)-\hat{x}^*(s;t_{k-1})}_Q \dopr{s}\cdot \norm{\tilde{x}(s;t_k)-\hat{x}^*(s;t_{k-1})}_Q^{\infty}+2\int_{t_{k}}^{t_{k-1}+T} \norm{\tilde{x}(s;t_k)-\hat{x}^*(s;t_{k-1})}_Q\dopr{s} \cdot\norm{\hat{x}^*(s;t_{k-1})}_Q^{\infty}$.
	\begin{equation} \label{eq:deltaJ1}
	\begin{split}
		\Delta \tilde{J}_1 \leq& \int_{t_{k}}^{t_{k-1}+T} \norm{\tilde{x}(s;t_k)-\hat{x}^*(s;t_{k-1})}_Q \dopr{s} \\ &\qquad\qquad\cdot \max_{s}{\norm{\tilde{x}(s;t_k)-\hat{x}^*(s;t_{k-1})}_Q} \\ &+2\int_{t_{k}}^{t_{k-1}+T} \norm{\tilde{x}(s;t_k)-\hat{x}^*(s;t_{k-1})}_Q\dopr{s} \\ &\qquad\qquad\qquad\cdot\max_{s}{\norm{\hat{x}^*(s;t_{k-1})}_Q}.
	\end{split}
	\end{equation}
	Using the result in Corollary~\ref{cor:cor1} and the robustness constraint~(\ref{eqInAETDRHC:robustconstraint}), it can be calculated that 
	%$A\leq \frac{\overline{\lambda}(Q)}{\underline{\lambda}(P)}\frac{L^2\beta T(1-\beta)T}{L\beta T-1}\delta(\frac{L^2\beta T}{L\beta T-1}\delta + 2[(1-\beta)M+\beta]\alpha\epsilon)$.
	\begin{equation}
		\begin{split}
			\Delta &\tilde{J}_1\leq \frac{\overline{\lambda}(Q)}{\underline{\lambda}(P)}\frac{L^2\beta T(1-\beta)T}{L\beta T-1}e^{L(1-\beta)T}\delta \\
			&\cdot \left[\frac{L^2\beta T}{L\beta T-1}e^{L(1-\beta)T}\delta + 2[(1-\beta)M+\beta]\alpha\epsilon\right].
		\end{split}
	\end{equation}
	For the second term 
	%$B = \int_{t_{k-1}}^{t_k}\normsq{\hat{x}^*(s;t_{k-1})}_Q + \normsq{\hat{u}^*(s;t_{k-1})}_R\dopr{s}$,
	\begin{equation}
		\Delta \tilde{J}_2 = -\int_{t_{k-1}}^{t_k}\normsq{\hat{x}^*(s;t_{k-1})}_Q + \normsq{\hat{u}^*(s;t_{k-1})}_R\dopr{s},
	\end{equation}
	it follows that 
	%$B \geq \int_{t_{k-1}}^{t_k}\normsq{\hat{x}^*(s;t_{k-1})-x(s;t_{k-1})+x(s;t_{k-1})}_Q\dopr{s}$.
	\begin{equation} \label{eqInTheoStability:DeltaJ2}
		\Delta \tilde{J}_2 \leq -\int_{t_{k-1}}^{t_k}\normsq{\hat{x}^*(s;t_{k-1})}_Q\dopr{s}.
	\end{equation}
	Since $x(t_0)\in \mathcal{X}\setminus\Omega{(\alpha\epsilon)}$, we can have $\Delta \tilde{J}_2 \leq -\frac{\underline{\lambda}(Q)}{\overline{\lambda}(P)}\beta T (\alpha\epsilon - \frac{L^2\beta T}{L\beta T-1}\delta)^2$. By Corollary~\ref{cor:cor1}, one can get the following result for the third term:
	\begin{equation} \label{eqInTheoStability:DeltaJ3}
		\Delta \tilde{J}_3 \leq \frac{L^4\beta^2 T^2}{(L\beta T-1)^2} e^{2L(1-\beta)T}\delta^2.
	\end{equation} 
	Consequently, it can be obtained that 
   	%$\Delta \tilde{J}\paren{x(s;t_k),u(s;t_k)}\leq A -B \leq -\frac{\underline{\lambda}(Q)}{\overline{\lambda}(P)(n+1)}\beta T (\epsilon - \delta)^2$ 
	\begin{equation}
		\begin{split}
			\Delta& \tilde{J}\paren{x(s;t_k),u(s;t_k)}\leq \Delta \tilde{J}_1 + \Delta \tilde{J}_2 + \Delta \tilde{J}_3 \\ &\leq -\frac{\underline{\lambda}(Q)}{\overline{\lambda}(P)(n+1)}\beta T (\alpha\epsilon-\frac{L^2\beta T}{L\beta T-1}\delta)^2
		\end{split}
   	\end{equation}
   	if the stability condition~(\ref{cond:stability1}) is satisfied.
%    \begin{equation}
%    \begin{split} 
%	    (L+\frac{1}{\beta T})\delta^2 + 2[(1-\beta)M+\beta]\alpha\epsilon\delta \leq \frac{n}{n+1}\beta T (\epsilon - (L+\frac{1}{\beta T})\delta)^2	
%	\end{split}
%    \end{equation}
	Due to the sub-optimality of the designed control $\tilde{u}(s;t_k)$ at $t_k$, we can achieve that the deceasing properties of the optimal cost function at $t_{k-1}$ and $t_k$ is guaranteed by 
	$
	  \Delta J\paren{\hat{x}^{*}(s;t_{k-1}),\hat{u}^{*}(s;t_{k-1})} \leq -\frac{\underline{\lambda}(Q)}{\overline{\lambda}(P)(n+1)}\beta T (\alpha\epsilon - \frac{L^2\beta T}{L\beta T-1}\delta)^2
	$,
	which consequently shows that the optimal cost functional $J^{*}$ is decreasing as $t$ approaches to infinity. Since the nominal state $\tilde{x}$ stays outside $\Omega{(\alpha\epsilon)}$, it thus follows that the lower bound for the decreasing of optimal functional $J^{*}$ is a positive constant. Assume that the nominal state $\tilde{x}$ cannot converge to the terminal set $\Omega{(\alpha\epsilon)}$ in finite time, then the optimal functional will decrease to $-\infty$ as time evolves to infinity, which is a contradiction to the fact that the optimal functional is quadratic. Similar argument can be found in~\cite{chen1998quasi,li2014event}.
	\par 
	Step 2: For all initial state $x(t_0)\in\Omega{(\alpha\epsilon)}$, we need to prove that the closed-loop system~(\ref{eq:closedloopSystem}) converges to $\Omega{(\bar{\epsilon})}$. 
%	By using the result of Lemma 1 in~\cite{chen1998quasi}, the derivative of Lyapunov function can be denoted as
%	\begin{equation*}
%	  \begin{split} 
%	  \dot{V}(x(t))&\leq -\normsq{x(t)}_{Q^*}+2x^{\top}(t)P\omega(t) \\
%	  &\leq \frac{\underline{\lambda}(Q^*)}{\overline{\lambda}(P)}(-\normsq{x(t)}_{P}+\bar{\epsilon}^2),
%	  \end{split} 
%	\end{equation*}
%	where $\bar{\epsilon}=\sqrt{\frac{2\overline{\lambda}(P)\norm{\sqrt{P}}}{\underline{\lambda}(Q^*)}\epsilon\rho}$. Thus the system is stable. 
	Using the fact that $\hat{x}^*(s;t_{k-1})\in \Omega(\alpha\epsilon)$, $\Delta \tilde{J}_1$ in~(\ref{eq:deltaJ1}) can be rewritten as
	\begin{equation} \label{eqInStability:deltaJ1small}
		\Delta \tilde{J}_1 \leq \frac{\overline{\lambda}(Q)}{\underline{\lambda}(P)}\frac{4 \alpha \epsilon L^2\beta T(1-\beta)T}{L\beta T-1} e^{L(1-\beta)T}\delta.
	\end{equation}
	For the second term $\Delta \tilde{J}_2$ in~(\ref{eqInTheoStability:DeltaJ2}), we have the following inequality hold by Lemma~\ref{lem:chen1998}:
	\begin{equation} \label{eqInStability:deltaJ2small}
		\Delta \tilde{J}_2 \leq - \frac{\underline{\lambda}(Q)}{\overline{\lambda}(P)} \beta T \normsq{\hat{x}^*(t_k;t_{k-1}) }_P.
	\end{equation}
	According to the event-triggering condition, we can have $\normsq{x(t_k)-\hat{x}^*(t_k;t_{k-1}) }_P \leq \left(\frac{L^2\beta T}{L\beta T-1}\delta \right)^2$ by following a similar procedure in Corollary~\ref{cor:cor1}. Then applying the above inequality to~(\ref{eqInStability:deltaJ2small}) yields 
	\begin{equation} \label{eqInStability:deltaJ2small2}
		\begin{split}
			\Delta \tilde{J}_2 &\leq - \frac{\underline{\lambda}(Q)}{\overline{\lambda}(P)} \beta T \left[\normsq{x(t_k)}_P - \left(\frac{L^2\beta T}{L\beta T-1}\delta \right)^2\right] \\
			&\leq -\frac{\underline{\lambda}(Q)}{\overline{\lambda}(P)}\beta T \normsq{x(t_k)}_{P} + \frac{\underline{\lambda}(Q)}{\overline{\lambda}(P)} \frac{L^4\beta^2 T^2}{(L\beta T-1)^2}\delta^2.
		\end{split}
	\end{equation}
	Combining~(\ref{eqInStability:deltaJ1small}),~(\ref{eqInStability:deltaJ2small2}) and~(\ref{eqInTheoStability:DeltaJ3}), it follows that $\Delta J \leq \Delta \tilde{J}_1 + \Delta \tilde{J}_2 + \Delta \tilde{J}_3 \leq -\frac{\underline{\lambda}(Q)}{\overline{\lambda}(P)}\beta T \normsq{x(t_k)}_{P} + \frac{\underline{\lambda}(Q)}{\overline{\lambda}(P)} \frac{L^4\beta^2 T^2}{(L\beta T-1)^2}\delta^2 + \frac{\overline{\lambda}(Q)}{\underline{\lambda}(P)}\frac{4 \alpha \epsilon L^2\beta T(1-\beta)T}{L\beta T-1} e^{L(1-\beta)T}\delta + \frac{L^4\beta^2 T^2}{(L\beta T-1)^2} e^{2L(1-\beta)T}\delta^2$, which implies that the sate will converge to the set $\Omega(\bar{\epsilon})$ with $\bar{\epsilon} =   \left( 1+\frac{\overline{\lambda}(P)}{\underline{\lambda}(Q)} e^{2L(1-\beta)T} \right)\frac{L^4\beta T}{(L\beta T-1)^2}\delta^2 +  \frac{\overline{\lambda}(P)}{\underline{\lambda}(Q)} \frac{\overline{\lambda}(Q)}{\underline{\lambda}(P)}\frac{4 \alpha \epsilon L^2 (1-\beta)T}{L\beta T-1} e^{L(1-\beta)T}\delta$. Then the proof is completed by summarizing Step 1 and Step 2.
\end{proof}
\else
\begin{proof}
	Due to the space limit, please refer to the detailed technical proof in our full version~\cite{sun2018integral}.
\end{proof}
\fi
\begin{remark}
	The inequality~(\ref{cond:stability1}) shows that the stability can be guaranteed by properly designing the prediction horizon $T$, the triggering level $\delta$, and the contraction rate $M$ for the robustness constraint. The larger triggering level $\delta$ leads to less frequent sampling events by sacrificing the control performance, whereas the larger prediction horizon $T$ usually provides better control performance due to the fact that longer state evolution is considered in the optimization. However, solving the MPC problem with larger prediction horizon $T$ consumes more computation resources. Thus, a trade-off should be considered when designing the parameters $T$ and $\delta$ in terms of the control performance.
\end{remark}

\section{Simulation Results} \label{sec:simulationresult}
Consider a nonlinear cart-damper-spring system with the following dynamics:
\begin{equation*} \label{eq:cartdamperspring}
\left\{
\begin{aligned}
&\dot{x}_1(t)=x_2(t),  \\
&\dot{x}_2(t)=-\frac{\tau}{M_c}e^{-x_1(t)}x_1(t)-\frac{h_d}{M_c}x_2(t)+\frac{u(t)}{M_c}+\frac{\omega(t)}{M_c},
\end{aligned}
\right.
\end{equation*}
where $x_1(t)$ denotes the displacement of the cart, $x_2(t)$ is the velocity, its mass $M_c=1.25\,\text{kg}$, the nonlinear factor $\tau=0.9\,\text{N/m}$, the damping factor $h_d=0.42\,\text{N*s/m}$, and the constrained control input $u(t)\in[-1,1]$. The Lipschitz constant $L$ is $1.4$. For this integral-type ET-MPC, we choose the weighted matrices
\ifFullVersion
\begin{figure}[!htbp]
	\centering
	\includegraphics[width=\linewidth]{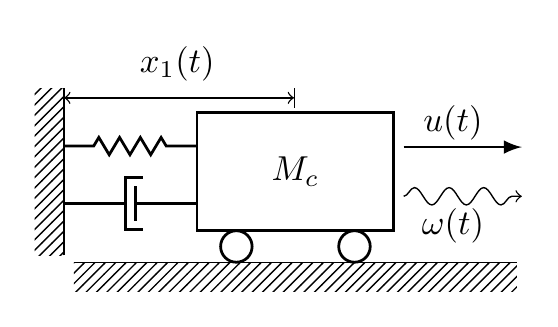}%
	\caption{The schematic illustration of a cart-damper-spring system.}
	\label{fig:cart}
\end{figure}
\fi
$Q=[0.1, 0.0; 0.0, 0.1]$ and $R=0.1$. Then the Linear Quadratic Regulator (LQR) feedback gain for~(\ref{eq:linearizedDynamics}) can be calculated as $K=[-0.4454,-1.0932]$. According to Lemma~\ref{lem:chen1998}, the corresponding $P$ matrix is designed as $P=[0.1692, 0.0572; 0.0572, 0.1391]$ and the terminal set level is determined as $\epsilon=0.03$. We choose the scaling ratio of the terminal set as $\alpha=0.8$ and the parameter as $\beta=0.6$. In addition, $T=2.0\,\text{s}$ and $M=10$ are chosen to satisfy the feasibility and stability conditions~(\ref{eqInTheo2:c1}),~(\ref{eqInTheo2:c2}),~(\ref{cond:stability1}). Therefore, the minimum inter-execution time is $\beta T=1.2\,\text{s}$. By using~(\ref{eqInTheo2:maxdisturbance}) in Theorem~\ref{theo:feasibilityResult}, the maximum allowable disturbance is calculated as $\rho_{\max}=0.00058$. Thus the additive disturbance can be configured as $\rho=0.00031$. The triggering level is chosen as $\delta=8.1\times 10^{-5}$ in order to satisfy the stability condition~(\ref{cond:stability1}).
% \par
% \thirdRevision{%
% Besides, the event-triggered control performance under different disturbance levels is investigated by considering the control performance and the communication performance, i.e., the finite horizon cost ($\int_{0}^{12}$$\normsq{x(s)}_{Q}+\normsq{u(s)}_{R}\dopr{s}$) and the average sampling numbers $\frac{12}{\text{Average Triggering Interval}}$). for different triggering level $\delta$ and prediction horizon $T$ is shown in Table~\ref{tab:tradeoff_triggeringlevel}. It can be seen from Table~\ref{tab:tradeoff_triggeringlevel} that the performance index decreases when the prediction horizon $T$ increases and increases when the triggering level $\delta$ decreases, which shows that a trade-off may exist for designing $\rho$.}
\par
\begin{figure}
	\centering
	\includegraphics[width=\linewidth]{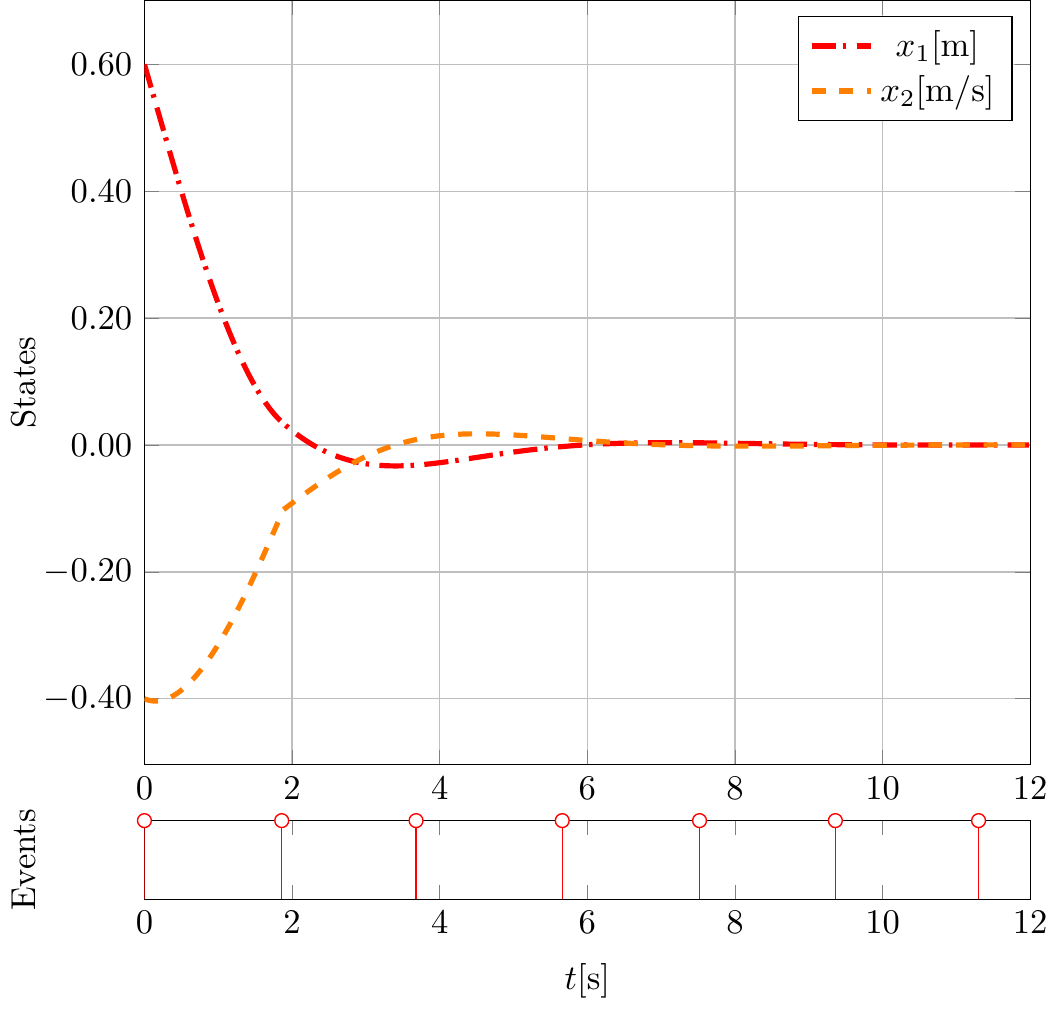}%
	\caption{States trajectories of the closed-loop system~(\ref{eq:closedloopSystem}) driven by integral-type ET-MPC~(\ref{eq:AETDRHC}), and event-triggering instants with condition~(\ref{eq:firstTrigger}). The red circle denotes the event-triggering instant.}
	\label{fig:system_trajectory}
\end{figure}
\begin{figure}
	\centering
	\includegraphics[width=\linewidth]{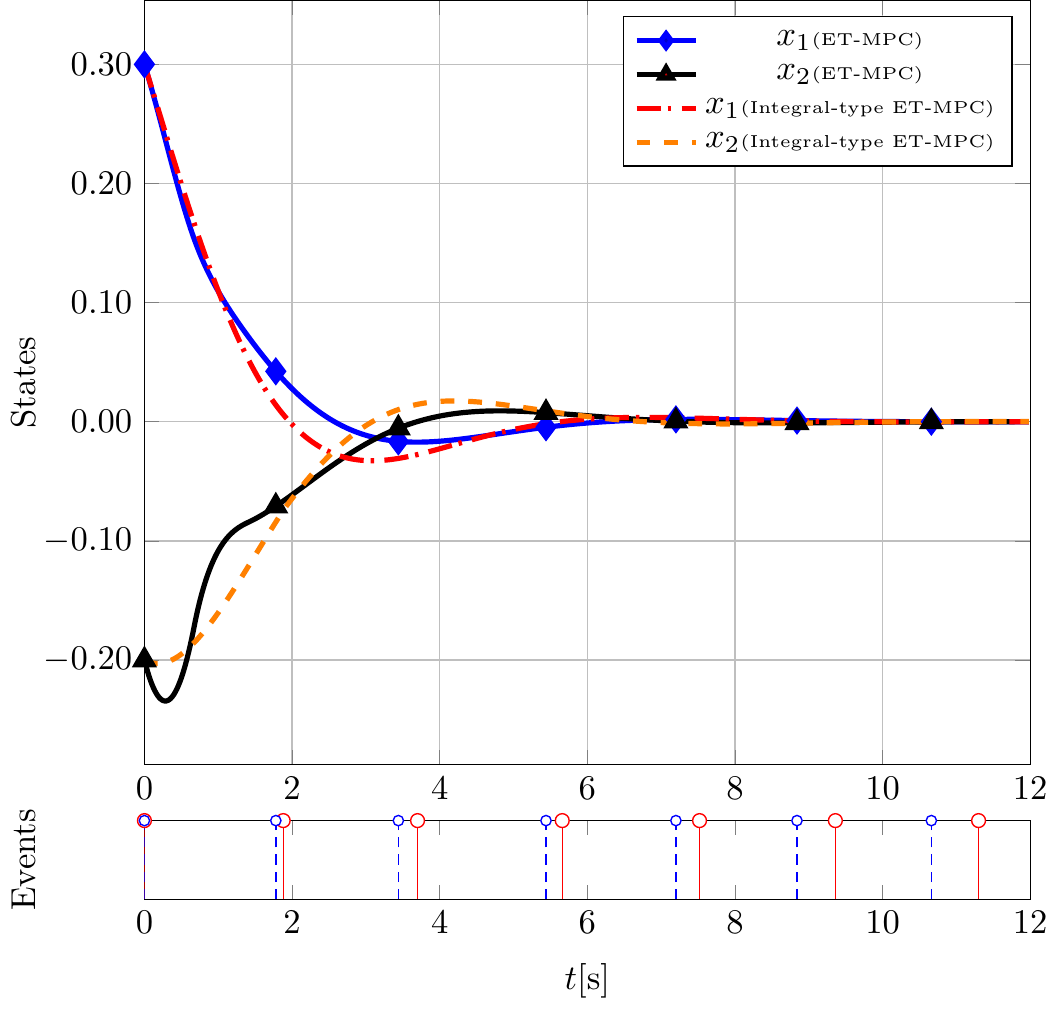}
	\caption{Comparison of state trajectories and event-triggering instants by the integral-type ET-MPC and the conventional ET-MPC in~\cite{li2014event}. The red circle denotes the event-triggering instant by our method and the blue one denotes the event-triggering instant by the conventional ET-MPC.}
	\label{fig:stateeventscomparison}
\end{figure}
We use IPOPT~\cite{wachter2006implementation} to solve the online optimization problem. Fig.~\ref{fig:system_trajectory} shows the state trajectory and event-triggering instants by using Algorithm~\ref{alg:int_ETRHC} given the initial state $x_0=[0.6, -0.4]$. Under the same initial state, the ET-MPC algorithm in~\cite{li2014event} is not feasible, which might indicate that the proposed integral-type ET-MPC scheme with less-conservative robustness constraint can admit an enlarged initial feasible region. For comparison purposes, we conduct another numerical example in Fig.~\ref{fig:stateeventscomparison}, where the initial states are set as $x_0=[0.3, -0.2]$ to satisfy the initial feasibility of both the conventional ET-MPC and the integral-type ET-MPC. In order to further show the advantages of the proposed integral-type ET-MPC scheme, we have also done two Monte-Carlo simulations to compare the communication performance with that of the conventional ET-MPC. The simulation results show that the average sampling frequency of our approach is $6.38$ which outperforms $6.66$ of the conventional ET-MPC. One can see that the integral-type ET-MPC can save considerable communication resource by performing less frequent event-triggered samplings, which is more efficient than the conventional ET-MPC.
% \newcommand{\tabincell}[2]{\begin{tabular}{@{}#1@{}}#2\end{tabular}}
% \begin{table}[htbp]
% 	\centering
% 	\caption{The event-triggered control performance of the integral-type ET-MPC under different disturbance levels.}
% 	\begin{tabular}{ccc}
% 		\toprule
% 		\tabincell{c}{Disturbance\\Levels} &  \tabincell{c}{Control\\Performance}  &  \tabincell{c}{Communication\\Performance} \\
% 		\midrule
% 		$0.00031$ &  $0.0132$ & $6.26$   \\
% 		$0.00062$ &  $0.0131$ & $9.02$   \\
% 		$0.00093$ &  $0.0131$ & $11.03$  \\
% 		$0.00124$ &  $0.0131$ & $12.61$	\\
% 		\bottomrule
% 	\end{tabular}
% 	\label{tab:tradeoff_triggeringlevel}
% \end{table}

\section{Conclusion and Future Work} \label{sec:conclusion}
In this paper, we have designed an integral-type ET-MPC scheme for nonlinear systems with additive disturbance. The integral-type ETM has shown considerable improvement on avoiding unnecessary communication. A new less conservative robustness constraint was proposed to handle the additive disturbance. For the feasibility and stability of the proposed integral-type ET-MPC framework, we have established the conditions for guaranteeing these two properties. In addition, we showed that the feasibility and stability properties are related to the prediction horizon, the disturbance bound, the triggering level, and the contraction rate for the robustness constraint. The simulation and comparison study demonstrated the effectiveness of proposed method. Future study will be focused on the output feedback ET-MPC and distributed MPC with the integral-type event-triggering mechanism.

\bibliographystyle{IEEEtran}
\bibliography{sqbib}

\end{document}